\documentclass[final]{siamltex}

\usepackage{graphicx}
\usepackage{amssymb,amsmath}
\usepackage{slashbox}
\graphicspath{ {c:/Users/Ronald_Morgan/PCTeXv4/FIGURES/} }

\begin{document}

\bibliographystyle{plain}

\title{Two-Grid Deflated Krylov Methods for Linear Equations\footnotemark[1]}
%Change GMRES to Krylov if defl.BiCGST works well.

\author{
Ronald B. Morgan\footnotemark[2],
Travis Whyte\footnotemark[3],
Walter Wilcox\footnotemark[4]
\and Zhao Yang\footnotemark[5]}

\maketitle

\renewcommand{\thefootnote}{\fnsymbol{footnote}}
\footnotetext[1]{The first author was supported by the National Science Foundation under grant DMS-1418677.  We also acknowledge support from a current Baylor University Research Committee Grant.}
\footnotetext[2]{Department of Mathematics, Baylor
University, Waco, TX 76798-7328 ({\tt Ronald\_Morgan@baylor.edu}).}
\footnotetext[3]{Department of Physics, Baylor
University, Waco, TX 76798-7316 ({\tt Travis\_Whyte@baylor.edu}).}
\footnotetext[4]{Department of Physics, Baylor
University, Waco, TX 76798-7316 ({\tt Walter\_Wilcox@baylor.edu}).}
\footnotetext[5]{Department of Mathematics, Baylor
University, Waco, TX 76798-7328 ({\tt zhao\_yang@alumni.baylor.edu}).}

\begin{abstract}
An approach is given for solving large linear systems that combines Krylov methods with use of two different grid levels.  Eigenvectors are computed on the coarse grid and used to deflate eigenvalues on the fine grid.  GMRES-type methods are first used on both the coarse and fine grids.  Then another approach is given that has a restarted BiCGStab  (or IDR) method on the fine grid.  While BiCGStab is generally considered to be a non-restarted method, it works well in this context with deflating and restarting.  Tests show this new approach can be very efficient for difficult linear equations problems.
\end{abstract}

%Note to Editors: This paper is about a new approach to solving large systems of linear equations that combines Krylov methods with some of multigrid.  It is novel in the way these are combined in order to deflate eigenvalues.  It also is unique in effectively restarting the normally non-restarted method BiCGStab.  Applications are given to several areas including QCD.

\begin{keywords}
 linear equations, deflation, GMRES, BiCGStab, eigenvalues, two-grid
\end{keywords}

\begin{AMS}
65F10, 15A06
\end{AMS}

\pagestyle{myheadings}
\thispagestyle{plain}
\markboth{R. B. MORGAN and Z. YANG}{Two-grid Deflated Krylov Methods}

\section{Introduction}

We look at solving large systems of linear equations $Ax = b$ that come from discretization of partial differential equations.  There are a variety of iterative methods for solving these problems.  In particular, multigrid methods~\cite{Fe64, Br77, BrennerScott, MultigrTut,TrOoSc} are extremely effective under certain conditions.  However, there are many situations for which standard multigrid methods do not work. Here we give a different approach to using the power of different grid levels for such situations.  This will be done by combining with Krylov subspace methods.

Convergence of iterative methods is generally affected by the conditioning of the linear equations, or more specifically by the presence of small eigenvalues.  One reason that multigrid methods are effective is that eigenvectors corresponding to small eigenvalues generally have similar shape on different grids levels, and multigrid is able to handle the small eigenvalues on coarse grids where both the linear equations problem is better conditioned and the computation is cheaper due to the smaller matrix.

Krylov subspace methods are more robust than multigrid in the sense that they can be applied to problems for which multigrid fails.  However, they can converge slowly.  Work has been done on dealing with the detrimental presence of small eigenvalues for Krylov methods.  Restarted methods such as GMRES are particularly sensitive to the presence of small eigenvalues.  Deflated GMRES methods~\cite{GMRES-E,KhYe,ChSa,BaCaGoRe,BuEr,FrVu,GMRES-DR,PadeStMaJoMa,GiGrPiVa,GaGuLiNa} compute approximate eigenvectors and use them to remove or deflate the effect of small eigenvalues.  In particular, we will use the method GMRES-DR~\cite{GMRES-DR} which both solves linear equations and simultaneously computes eigenvectors.

Once approximate eigenvectors have already been computed, one deflation technique is to use them to build a preconditioner for the linear equations~\cite{KhYe,BaCaGoRe,BuEr,FrVu,PadeStMaJoMa,SiEmMo}.  However, such a preconditioner can become expensive if there are many eigenvalues being deflated; each eigenvector is applied at every iteration.  Deflation can be done with less expense outside of the Krylov iteration~\cite{gproj,StOr}.  We will use such a method called GMRES-Proj~\cite{gproj} which projects over approximate eigenvectors in between cycles of regular GMRES.

The computation of approximate eigenvectors can be expensive for difficult problems.
The new approach in this paper is to compute eigenvectors on a coarse grid and move them to the fine grid.  GMRES-DR is used on the coarse grid to both find approximate eigenvectors and generate an approximate linear equations solution.  Then GMRES-Proj takes what has been found by GMRES-DR and solves the fine grid linear equations.

For fairly sparse matrices, the orthogonalization expense in GMRES-Proj can be a major expense.  We present an approach that substitutes either BiCGStab or IDR in place of GMRES.  Here we restart BiCGStab and IDR even though they are normally non-restarted methods.  This may seem risky, because of the inconsistent convergence of these nonsymmetric Lanczos methods.  However, in our experiments, restarted BiCGStab and IDR converge reliably.

For related work, see Elman, Ernst and O'Leary~\cite{ElErOL} for another way of combining multigrid and GMRES.
Closer to this paper is Sifuentes' thesis~\cite{Si10} which has two-grid deflation, but with Arnoldi on the coarse grid and with the more expensive approach of building a deflating preconditioner for GMRES on the fine grid.
Instead of using computed approximate eigenvectors, Erlangga and Nabben~\cite{ErNa08A} deflate using vectors from the interpolation operator that maps from coarse to fine grid.

Section 2 reviews some of the previous methods that will be used. Section 3 presents the method Two-grid Deflated GMRES.  Two-grid Deflated BiCGStab/IDR is then given in Section 4.  Further examples are in Section 5 with a Helmholtz example that includes multigrid preconditioning and with an example of improving eigenvectors on the fine grid.

\section{Review}

Here we very quickly describe methods that will be used in the rest of the paper.

\subsection{GMRES-DR}

The GMRES with Deflated Restarting (GMRES-DR)~\cite{GMRES-DR} method uses Krylov subspaces to both solve linear equations and compute the eigenpairs with smallest eigenvalues.  Once eigenvectors converge far enough, their presence in the subspace can essentially remove or deflate the effect of the small eigenvalues on the linear equations.  Usually only moderate accuracy is needed before the approximate eigenvectors have a beneficial effect.

GMRES-DR(m,k) saves $k$ approximate eigenvectors at the restart and builds out to a subspace of dimension $m$.  Specifically, for one restarted cycle it uses subspace
\[Span\{y_1, \ldots, y_k, r_0, Ar_0, A^2 r_0, \ldots A^{m-k-1} r_0 \},\]
where $y_i$'s are harmonic Ritz vectors from the previous cycle and $r_0$ is the residual vector at the start of the cycle.  This augmented subspace is actually a Krylov subspace itself, and it contains Krylov subspaces with each $y_i$ as starting vector.  This makes the eigenvectors generally converge along with the linear equations.

GMRES-DR converges faster than restarted GMRES for difficult problems with small eigenvalues.  It also often converges faster than BiCGStab in terms of matrix-vector products, but has greater orthogonalization costs per matrix-vector product.

\subsection{GMRES-Proj}

There are situations where approximate eigenvectors are available at the beginning of the solution of linear equations.  For example, if there are multiple right-hand sides, then eigenvectors could have been computed during solution of earlier right-hand sides~\cite{PadeStMaJoMa,gproj,StOr}.
The method GMRES-Proj~\cite{gproj} uses these approximate eigenvectors to deflate the corresponding eigenvalues while solving linear equations.
GMRES(m)-Proj(k) assumes that $k$ approximate eigenvectors have been previously computed and alternates projection over these vectors with cycles of GMRES(m).

\vspace{.10in}
\begin{center}
\textbf{GMRES(m)-Proj(k)}
\end{center}
\begin{description}
 \item[0.] Let $k$ be the number of approximate eigenvectors that are available.  Choose $m$, the dimension of subspaces generated by restarted GMRES.
 \item[1.] Alternate between A) and B) until convergence:
 \begin{description}
 \item[A)] Apply Galerkin projection over the subspace spanned by the $k$ approximate eigenvectors.
 \item[B)] Apply one cycle of GMRES(m).
\end{description}
\end{description}
\vspace{.15in}

The algorithm for the projection step is given next.

\vspace{.10in}
\begin{center}
\textbf{Galerkin projection over a set of approximate eigenvectors}
\end{center}
\begin{description}
 \item[0.] Let the current system of linear equations be $A(x-x_0) = r_0$.
 \item[1.] Let $V$ be an $n$ by $k$ orthonormal matrix whos columns span the set of approximate eigenvectors.
 \item[2.] Form $H = V^T A V$ and $c = V^T r_0$.
 \item[3.] Solve $H d = c$, and let $\hat x = V d$.
 \item[4.] New approximate solution is $x_p = x_0 + \hat x$, and new residual is $r = r_0 - A \hat x = r_0 - A V d$.
\end{description}
\vspace{.15in}

MinRes projection can be used instead of Galerkin.  It is the same as Galerkin, except $H = (AV)^T A V$ and $c = (AV)^T r_0$.  Often Galerkin is best and can even be much better, however we have seen situations for which MinRes works better.  All experiments in this paper use Galerkin.

\subsection{Two-grid Arnoldi}

A two-grid method for computing eigenvalues and eigenvectors is given in~\cite{MGArn}.
Eigenvectors are computed on a coarse grid with a standard Arnoldi method, are moved to the fine grid (with spline interpolation), then are improved on the fine grid with Arnoldi-E~\cite{Arnoldi-R}, a method that can accept initial approximations.

\section{Two-grid Deflated GMRES}

Our new methods generate approximate eigenvectors from the coarse grid and use them to deflate eigenvalues on the fine grid.  In this section, we give a version using GMRES methods.  GMRES-DR is applied on the coarse grid.  This generates eigenvectors, and solves the coarse grid linear equations.  This linear equations solution is moved to the fine grid with spline interpolation or prolongation and used there as the initial guess.  The eigenvectors are similarly moved to the fine grid and, if necessary, improved on the fine grid.  GMRES-Proj is applied on the fine grid using these approximate eigenvectors.
This can have much faster convergence than restarted GMRES.  Compared to running GMRES-DR on the fine grid, it is cheaper to implement and can deflate eigenvalues from the beginning.

\vspace{.10in}
\begin{center}
\textbf{Two-grid Deflated GMRES}
\end{center}
\begin{description}
 \item[0.]  Choose $m$ and $k$ for the coarse grid.  Pick $nev$, the number of eigenpairs that are required to converge to an eigenvalues tolerance, say $rtolev$.  For the fine grid, pick the linear equations residual tolerance, $rtol$.  Also choose $m3$ for GMRES on the fine grid.
 \item[1.] Apply GMRES-DR(m,k) on the coarse grid.  Move the approximate eigenvectors to the fine grid (we generally use spline interpolation).  Move the solution of the coarse grid linear equations problem to the fine grid and use it as an initial guess for the fine grid problem.
 \item[2.] (If needed:) Improve approximate eigenvectors on the fine grid using Arnoldi-E (see Two-grid Arnoldi method in~\cite{MGArn}).
 \item[3.] Apply GMRES(m3)-Proj(k) on the fine grid.
\end{description}
\vspace{.15in}

For the first examples there is no need for the second phase, but it is used in Subsections 5.3 and 5.4.

{\it Example 1.} We consider a system of linear equations from finite difference discretization of the 2-D convection-diffusion equation $ - \operatorname{e}^{5xy} (u_{xx} + u_{yy}) + 40 u_{x} + 40 u_y = c \thinspace sin x \thinspace cos x \operatorname{e}^{xy}$ on the unit square with zero boundary conditions and $c$ chosen to make the right-hand side be norm one.  The discretization size is $h = \frac{1}{512}$, leading to a matrix of dimension $n = 511^2 = 262,121$.  The coarse grid discretization size is $h = \frac{1}{64}$ giving a matrix of dimension $63^2 = 3969.$

The first phase on the coarse grid uses GMRES-DR(150,100) and runs until $80$ eigenpairs have converged to a level of residual norm below $10^{-8}$.  These residual norms are computed only at the end of cycles.  The eigenvectors are moved to the fine grid and are accurate enough there to be effective in deflating eigenvalues (after the Rayleigh-Ritz procedure is applied to all 100 vectors that are moved from the coarse grid, the smallest 80 Ritz pairs have residual norms at or below $1.4*10^{-3}$).  So as mentioned, the second phase is not needed.  The third phase of solving the fine grid linear equations is stopped when the relative residual norm drops below tolerance of $rtol = 10^{-10}$.  This is also checked only at the end of each GMRES cycle, but of course could be easily monitored during the GMRES runs.

The top of Figure 3.1 has convergence of both linear equations and eigenvalues on the coarse grid.  The linear equations converge before any of the eigenvalues, and it actually takes much longer for all 80 eigenpairs to become accurate.  The linear equations converge in 19 cycles of GMRES-DR(150,100) which use 1050 matrix-vector products (150 for the first cycle and 50 each for the next 18).  The eigenvalues take 107 cycles or 5450 matrix-vector products.

The bottom of the figure shows convergence of the linear equations on both grids versus the number of fine-grid-equivalent matrix-vector products.  The coarse grid matrix is about 64 times smaller than the the fine grid matrix, so we scale the number of matrix-vector products by a factor of 64 to get the fine-grid-equivalents.  The coarse grid linear equations then converge so rapidly that the convergence curve is barely noticeable on the left of the graph.  Then there is a small gap before the curve for fine grid convergence starts.  This gap is for both the convergence of the eigenpairs on the coarse grid and the matrix-vector products needed to form the projection matrix $H$ for the fine grid.
Three different values of $m3$ are used for GMRES(m3)-Proj(100): 50, 100 and 200.
While the fastest convergence is with GMRES(200), the least expensive is for $m3 = 100$.   We define the approximate cost as $cost = 5*mvp + vops$, where $mvp$ is the number of matrix-vector products, the 5 comes from the approximate number of non-zeros per row, and $vops$ is the number of length-$n$ vector operations such as dot-products and daxpy's.    The costs for the entire process, including the first phase on the coarse grid, is $cost = 1.10*10^6$ for $m3 = 50$, $1.06*10^6$ for $m3=100$ and finally $1.30*10^6$ with $m3=200$. The orthogonalization cost is high for GMRES(200) and this is the biggest expensive for a sparse matrix such as this one.

\begin{figure}
\vspace{-2.65in}
\hspace{-.7in}
\includegraphics[scale=.75]{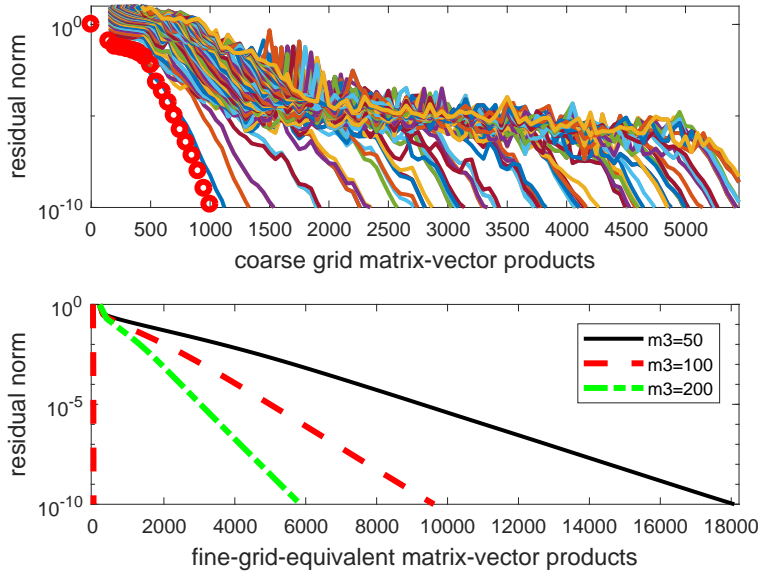}
\vspace{-2.8in}
\caption{This figure has Two-grid Deflated GMRES for the convection-diffusion example.  Fine grid matrix size is $n=261,121$ and coarse grid matrix size is $n=3969$.
The top portion has GMRES-DR(150,100) on the coarse grid; linear equations convergence is shown with circles at the end of each cycle and convergence of 80 eigenpairs is shown with lines.  The bottom portion has the coarse grid linear equations solution on the very left, scaled by 64 to correspond to fine-grid matrix-vector products.  Then this is followed by GMRES(m3)-Proj(100) on the fine grid, with $m3 = 50, 100, 200$. }
\end{figure}

Next, we consider different sizes for the coarse grid.  Smaller coarse grid matrices mean less work is needed to find the eigenpairs, however they may not be as accurate for the fine grid work.  Table 3.1 has coarse grids from size $65,025 = 255^2$ down to $49 = 7^2$.  All of the tests use the same type coarse grid computation as before and use GMRES(100)-Proj(100) on the fine grid.  The results show that for this matrix, the fine grid convergence is fairly robust with respect to the coarse grid size.  The number of fine grid cycles increases by less than half as the coarse grid goes from size 65,025 down to 255.  For the 255 case, the accuracy of 80 eigenpairs on the fine grid is $5.7*10^{-3}$ or better, and that is enough to be fairly effective.  The convergence is 10 times faster than with no deflation which is on the last row of the table.  The larger coarse grid matrices do give more accurate eigenvectors on the fine grid, for instance with residual norms $4.4*10^{-5}$ and below for size $65,025$.  However, this greater accuracy is not really needed and only speeds up the fine grid convergence from 94 to 85 cycles compared to size $3969$.

\begin{table}

\caption{Effect of Coarse Grid Size}

\begin{center}
\begin{tabular}{|c|c|c|c|c|c|c|c|c|}  \hline\hline
 coarse grid    & coarse grid   & accuracy of 80   & fine grid & cost      \\
 matrix size    & cycles        & eigenpairs on f.g.  & cycles & (millions)\\ \hline \hline
$65,025 = 255^2$ & 1023         & $4.4*10^{-5}$ & 85        & 8.02      \\ \hline
$16,129 = 127^2$ & 270          & $3.1*10^{-4}$ & 86        & 1.40      \\ \hline
 $3969 = 63^2$  & 107           & $1.4*10^{-3}$ & 94        & 1.06      \\ \hline
 $961 = 31^2$   & 43            & $5.8*10^{-3}$ & 106       & 1.15      \\ \hline
 $225 = 15^2$   & 12            & $5.7*10^{-3}$ & 120       & 1.30      \\ \hline
 $49 = 7^2$     & 1             & ($9.2*10^{-3}$ for 40) & 177  & 1.91  \\ \hline
 no coarse grid & -             & -             & 1255      & 13.2      \\ \hline
    \hline

\end{tabular}
\end{center}
%\label{}

\end{table}

Now we discuss other methods.  Multigrid is much faster than Krylov when convection is low.  However, here with convection terms  $40 u_{x} + 40 u_y$, standard multigrid methods do not converge.
Next we compare to other Krylov methods.
The top part of Figure 3.2 has convergence in terms of matrix-vector products.  GMRES-DR(150,100) is run just on the fine grid, as are IDR and BiCGStab.  The new Two-grid Deflated GMRES method is given with $m3 = 100$ and it converges faster than fine grid GMRES-DR(150,100), because approximate eigenvectors from the coarse grid allow deflation of small eigenvalues from the beginning.  The new method also converges in many less matrix-vector products than the non-restarted methods BiCGStab~\cite{vdV92} and IDR~\cite{SovG,SovG11}.  However, for this quite sparse matrix, the new method is not necessarily better overall than the non-restarted methods.  The bottom of Figure 3.2 has convergence as a function of the approximate cost and BiCGStab is much faster because of using only seven $vops$ per $mvp$.  IDR(4) uses 12.8 $vops$ per $mvp$, but is not effective for this matrix.

\begin{figure}
\vspace{-2.65in}
\hspace{-.7in}
\includegraphics[scale=.75]{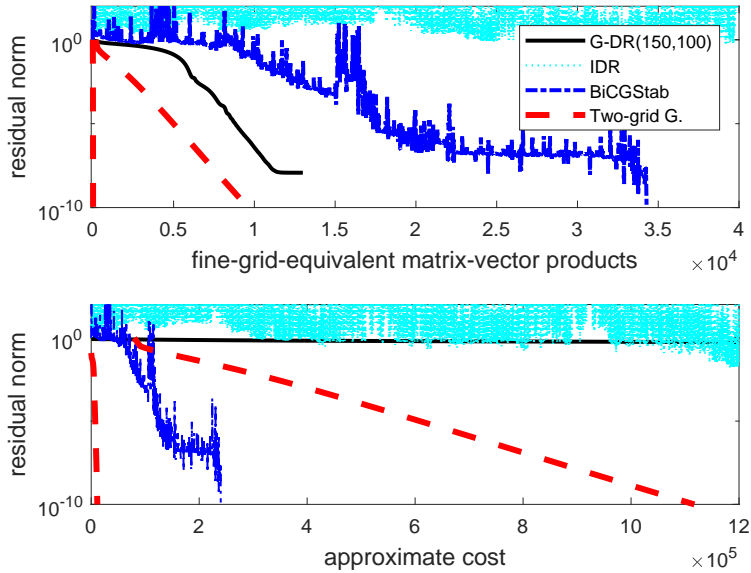}
\vspace{-2.8in}
\caption{Two-grid deflated GMRES (with GMRES-DR(150,100) on the coarse grid and GMRES(100)-Proj(100) on the fine grid) compared to GMRES-DR(150,100), BiCGStab and IDR(4). The top half is plotted against matrix-vector products and the bottom half against the cost in terms of vector op's.}
\end{figure}

We have seen in this example that Two-grid Deflated GMRES can be better than BiCGStab in terms of matrix-vector products.  So it can be the better method if the matrix is not very sparse or if an expensive preconditioner is used.  Then the number of iterations is the main factor instead of the GMRES orthogonalization expense.

However, for sparse problems, BiCGStab has an advantage.  This motivates the next method. For the fine-grid portion of the two-grid approach, GMRES is replaced by BiCGStab.  This is designed to have both a low number of matrix-vector products and low costs for vector operations.

\section{Two-grid Deflated BiCGStab}

\subsection{The algorithm}

We wish to use approxiate eigenvectors from the coarse grid to deflate eigenvalues from BiCGStab or IDR on the fine grid.  In~\cite{defbi,NLan-DR,AROrSt} a deflated BiCGStab method is given.  A single projection is applied before running BiCGStab using both right and left eigenvectors.  Here a single projection will not be effective, because our eigenvectors are not accurate on the fine grid.  So we implement BiCGStab/IDR as a restarted method with projections at each restart.  Not only does this allow us to use less accurate eigenvectors, but also it does not require left eigenvectors.  We use the same Galerkin projection as for Two-grid Deflated GMRES.

We now give the implementation of this restarted, deflated BiCGStab/IDR method.  It replaces phase 3 in the Two-grid Deflated GMRES algorithm given earlier and thus is part of a Two-grid Deflated BiCGStab/IDR algorithm.  The new BiCGStab(ncyc)-Proj(k) is similar to GMRES(m3)-Proj(k), but $ncyc$ is the total number of BiCGStab cycles, not the length of the cycles (so there are $ncyc-1$ restarts).  As before, the $k$ gives the number of approximate eigenvectors that are being projected over.  The stopping test for the algorithm uses the minimum of two different quantities, but this means the maximum convergence between the two criteria.  First, we require each cycle to converge to a fraction of the remaining distance to final convergence, in terms of orders of magnitude.  This fraction is determined by the number of remaining cycles to equalize the amount expected of each cycle.  Second, for cycle $icyc$ we want to reach a point that is at least $icyc/ncyc$ of the way toward convergence (again in orders of magnitude).  The first criteria usually asks for convergence to a further point and so is the one enforced, but in the case where the residual norm jumps up during the projection, the second criteria is needed (see Example 6).

\vspace{.10in}
\begin{center}
\textbf{BiCGStab(ncyc)-Proj(k) and IDR(ncyc)-Proj(k)}
\end{center}
\begin{description}
 \item[0.] Assume $k$ approximate eigenvectors are provided.  Let $rtol$ be the specified relative residual tolerance for the linear equations solution.  Choose $ncyc$, the requested number of cycles of BiCGStab/IDR.  Let $r_0$ be the initial residual for the fine grid iteration.
 \item[1.] Apply GMRES-DR(m,k) on the coarse grid.  Move approximate eigenvectors to the fine grid and also move the solution of the coarse grid linear equations as an initial guess for the fine grid problem.
 \item[2.] (If needed:) Improve eigenvectors on fine grid using Arnoldi-E (see~\cite{MGArn}).
 \item[3.] For $icyc = 1:ncyc$
    \begin{description}
    \item[a)] Apply Galerkin projection over approximate eigenvectors.
    \item[b)] Let $\|r\|$ is the current residual.  Set the relative residual tolerance for this cycle, $rticyc$, to be the minimum (the further convergence point) of $(rtol*\|r_0\|/\|r\|)^{(\frac{1}{ncyc-icyc+1})}$ and $(\|r_0\|/\|r\|)*(rtol)^{\frac{icyc}{ncyc}}$.
    \item[c)] Run BiCGStab or IDR with relative residual tolerance of $rticyc$.
    \item[d)] Break out of the loop if $\|r\|$ is already below $rtol$.
    \end{description}

\end{description}
\vspace{.15in}

For our tests, the Matlab BiCGStab program is called.  For IDR, we use the program described by van Gijzen and Sonneveld in~\cite{SovG11} and available in MATLAB code from the authors.  The default version IDR(4) is called.

{\it Example 2. }
We return to the same convection-diffusion problem as in Example 1 that has $n = 262,121$ and coarse grid matrix of size $3969$.  Figure 4.1 has plots of BiCGStab(ncyc)-Proj(100) with $ncyc = 5$, 10 and 20.  The more frequent restarts, and thus more deflations of the eigenvalues, allow $ncyc=20$ to converge faster.  It is surprising that in spite of the very jagged behavior of the residual norms of BiCGStab, the overall convergence with $ncyc=20$ is quite consistent.  With even more frequent restarts, the convergence is similar to $ncyc=20$.  These tests are not shown on the figure, because they mostly overlie the $ncyc=20$ curve.  However, see Table 4.1 for the number of matrix-vector products needed with some other values of $ncyc$.  Note that with large values of $ncyc$, there may be convergence before all cycles are used and the break out of the loop in part d) of the algorithm is activated.  This is shown in parentheses in the table.  For example, only 168 runs are needed when 200 are specified.
We also tried restarting BiCGStab after 272 matrix-vector products (the average length for the $ncyc=20$ test), and the convergence was a little slower than the $ncyc=20$ method, using 5865 matrix-vector products instead of 5421.  Probably the $ncyc=20$ test has an advantage because it often restarts after an iteration where the residual norm comes down significantly.

Figure 4.2 has a comparison of different methods.  The top portion has convergence with respect to matrix-vector products.  Two-grid Deflated BiCGStab uses $ncyc=20$ and has projection over 100 approximate eigenvectors.  These vectors come from the same run of GMRES-DR(150,100) on the coarse grid as is used for Two-grid GMRES.  This deflated BiCGStab converges faster than Two-grid Deflated GMRES.
Both deflated methods are much better than regular BiCGStab.
Also given is a Two-grid Deflated BiCGStab result for solving a second system of linear equations with the same matrix but a different right-hand side (this right-hand side is generated randomly).  It converges similarly to the first right-hand side, but does not need the coarse-grid work or forming of the projection matrix.
The bottom of Figure 4.2 has a plot of convergence versus cost.  Two-grid BiCGStab converges faster than regular BiCGStab in spite of getting a late start due to the cost of the coarse grid phase.  The second right-hand side does not have this cost and so has less than half the expense.

\begin{figure}
\vspace{-2.65in}
\hspace{-.7in}
\includegraphics[scale=.75]{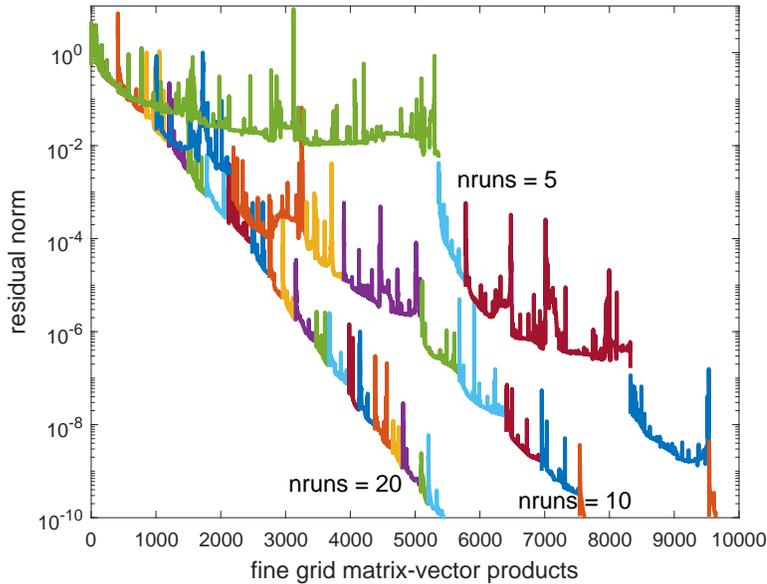}
\vspace{-2.8in}
\caption{Two-grid Deflated BiCGStab with GMRES-DR(150,100) on the coarse grid and BiCGStab(ncyc)-Proj(100) on the fine grid.  The number of cycles for the restarted BiCGStab phase is $ncyc = 5$, 10 and 20.  The color changes with each new cycle.}
\end{figure}

\begin{figure}
\vspace{-2.65in}
\hspace{-.7in}
\includegraphics[scale=.75]{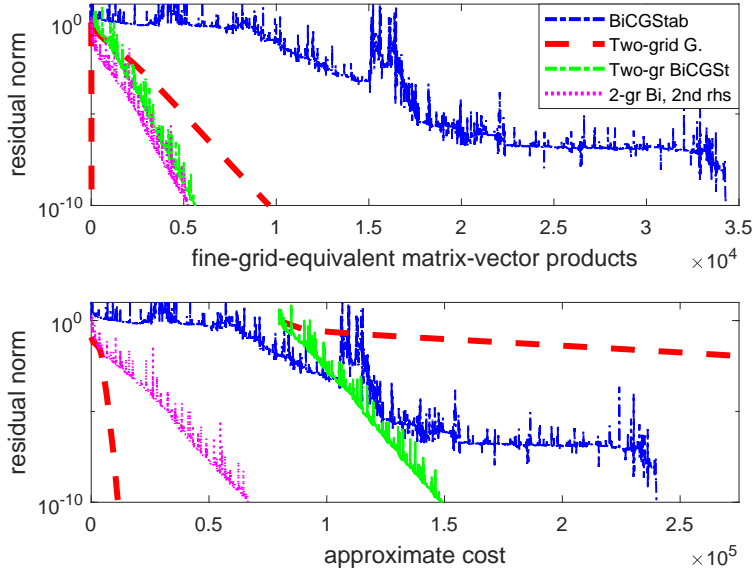}
\vspace{-2.8in}
\caption{Two-grid Deflated BiCGStab, with GMRES-DR(150,100) on the coarse grid and BiCGStab(20)-Proj(100) on the fine grid, compared to other methods.  The Two-grid Deflated GMRES has GMRES-DR(150,100) on the coarse grid and GMRES(100)-Proj(100) on the fine grid.  A test of deflated BiCGStab with a second right-hand side is also shown.}
\end{figure}

\begin{table}

\caption{Comparing different values of $ncyc$ for BiCGStab(ncyc)-Proj(100).  The number of matrix-vector products is given.}

\begin{center}
\begin{tabular}{|c|c|c|c|c|c|c|c|c|c|}  \hline\hline
 ncyc      & 5     & 10    & 15    & 20    & 30    & 50    & 100       & 150       & 200       \\
            &       &       &       &       &       &       &  (use 99) & (140) & (168) \\ \hline
 mvp's      &  9638 & 7606  & 6390  & 5421  & 5357  & 5278  & 5262  & 5591 & 5541     \\ \hline \hline

\end{tabular}
\end{center}
%\label{}

\end{table}

\subsection{Effectiveness of restarted BiCGStab}

%((NOTE: decided to not count extra mvp in tables, but am counting in the plots - funny, huh? ))

Figure 4.1 and Table 4.1 show something remarkable for deflated, restarted BiCGStab.  For $ncyc \ge 20$, the results are fairly close to being invariant of the number of cycles.  This is in spite of much smaller subspaces being used for the larger values of $ncyc$.  Also, the convergence is at a very consistent pace considering the usual erratic convergence of BiCGStab.
Deflated GMRES in Figure 3.1 is very sensitive to the subspace sizes, with $m3=50$ converging three times slower than for $m3=200$.  Deflated BiCGStab with 20 cycles uses an average of 271 matrix-vector products per cycle, while with $ncyc=100$, an average of 53 per cycle converges slightly faster.

We next investigate the reason for this effectiveness of restarted BiCGStab.  First, it is not due to the deflation.  To show this, the next example does not use deflation and has a similar phenomenon for restarted BiCGStab.

{\it Example 3.}  The matrix is the same as in the previous examples, except the size is $127^2 = 16,129$ and there is not a coarse grid matrix.  We use a random right-hand side to try to make the example more general, though generally there is not much effect from the right-hand side.  We run the restarted BiCGStab with no deflation between cycles.  Table 4.2 shows tests with different numbers of cycles.  Surprisingly, even for large numbers of cycles, the results are similar and often better than for regular non-restarted BiCGStab.  This behavior is reminiscent of Tortoise and the Hare behavior for GMRES~\cite{Em03} where smaller subspaces give faster convergence.  This is less surprising for BiCGStab, because of the lack of the minimal residual property.  However, in one way the behavior here is more extreme than Tortoise and the Hare for GMRES, because the convergence with frequent restarts can be even better than with no restarting.  For example, with $ncyc$ specified to be 400, the method ends up using 198 cycles, 3452 matrix-vector products and the average number of matrix-vector products per cycle of $17 \frac{1}{2}$.  This compares to 3969 matrix-vector products with non-restarted BiCGStab.  Also for GMRES(18), 28,838 matrix-vector products are required.

\begin{table}

\caption{Restarted BiCGStab with different numbers of cycles (no deflation is used).  The number of matrix-vector products is given.}

\begin{center}
\begin{tabular}{|c|c|c|c|c|c|c|c|c|c|}  \hline\hline
 ncyc   &  1 (no  & 5     & 10    & 25    & 50    & 75     & 150   & 250   & 400     \\
        & restart) & &  &       &       & (use 74) & (118) & (171) & (198)   \\ \hline
 mvp's  & 3969 & 3755  & 4371  & 3968  & 3386  & 3578   & 4049  & 3764  & 3452 \\ \hline \hline

\end{tabular}
\end{center}
%\label{}

\end{table}

It has been proposed that GMRES should use changing cycle lengths~\cite{BaJeKo}.  We suggest that variable cycle lengths is the main reason that restarted BiCGStab is so effective even with small subspace sizes.  Table 4.3 has results for restarted BiCGStab with average subspace dimensions of approximately 50, 35 and 18 (the lowest we could get was about 18, because when a high number of cycles is specified, the algorithm finishes in many fewer cycles than requested).  The number of matrix-vector products is similar for each of the three average subspace sizes (these three numbers come from Table 4.2 with $ncyc$ requested to be 75, 150 and 400).  The next column has GMRES(m) for $m = 50, 35$ and $18$.  As mentioned, the number of matrix-vector products is very large for small $m$.  The third result is for GMRES restarted exactly as BiCGStab is (so with the same number of matrix-vector products for each corresponding cycle, though it does not necessarily need all of the cycles).  Now GMRES uses even less matrix-vector products than BiCGStab.  This is a remarkable improvement over regular restarted GMRES.  Regular GMRES(m) can get stuck in a pattern~\cite{BaJeMa,CC} and changing the sizes of the subspaces can break the pattern.  The final column of the table has a different way of implementing variable restarting for GMRES.  A maximum cycle length is specified and the length of each individual cycle is randomly chosen between 1 and the max length (for example, max cycle length of 36 is used in the test because it gives an average of 18).  This approach performs much better than the fixed cycle length approach, but is not as good as using the BiCGStab lengths.  We conclude from these tests that our BiCGStab method restarts in a surprisingly effective way.
It is important for the new Two-grid Deflated BiCGStab method that frequent restarts still give an effective method, because if only a few restarts were used, there may not be frequent enough deflations, as seen in Figure 4.1 with the $ncyc = 5$ and $10$ cases.

Tables 4.2 and 4.3 make it seem like the sizes of the subspaces used for BiCGStab and GMRES, restarted as for BiCGStab, do not have much effect for this example.  This would be surprising, since it is well known that large subspaces can be an advantage for Krylov methods.  And indeed, even though GMRES does well with these small subspaces, it is significantly better without restarting, using 1383 matrix-vector products.  Since non-restarted BiCGStab has 3969 matrix-vector products, we see that here GMRES is better able to take advantage of a large non-restarted Krylov subspace.

\begin{table}

\caption{Matrix-vector products for restarting BiCGStab and GMRES.}

\begin{center}
\begin{tabular}{|c|c|c|c|c|}  \hline\hline
 average   & restarted     & GMRES(m),  &  GMRES restarted  & GMRES with \\
 cycle length         & BiCGStab      & m fixed   &  as was BiCGStab  & random restarts \\ \hline
 50     & 3578  & 6250      & 3266  & 3498      \\ \hline
 35     & 4049  & 11,760    & 3244  & 3733      \\ \hline
 18     & 3452  & 28,838    & 3411  & 5013      \\ \hline
 \hline

\end{tabular}
\end{center}
%\label{}

\end{table}

%((Zhao asks how restarted BiCGStab (without deflation) compares for Example 2 to non-rest Bi and deflated Bi. ))

\section{Further Experiments}

\subsection{Helmholtz matrix}
\

{\it Example 4.}  We next test the Helmholtz matrix from finite difference discretization of the differential equation $ - u_{xx} - u_{yy} - 100^2 u  = f$ on the unit square with zero boundary conditions.  The discretization size is again $h = \frac{1}{512}$ and again $n = 262,121$.  The $rtol$ is $10^{-10}$.  The right-hand side is generated random normal.
The coarse grid discretization size is $h = \frac{1}{128}$, giving a matrix of dimension $127^2 = 16,129.$
This problem is difficult because of the significantly indefinite spectrum.  Standard multigrid methods do not converge, however multigrid can be used as a preconditioner; see the next example.
Figure 5.1 has results for Krylov methods plotted against matrix-vector products.  We have found that IDR performs better than BiCGStab for indefinite matrices and only it is shown on the figure, but here even IDR does not converge.  GMRES-DR(200,150) converges eventually (though not quite to the requested $rtol$), however it is expensive due to orthogonalization costs.  Our Two-grid Deflated GMRES method first uses GMRES-DR(200,150) on the coarse grid until 120 eigenvalues converge to residual norm below $10^{-8}$.  This takes 285 cycles.  Then on the fine grid it has GMRES(100)-Proj(150).  This method converges faster than GMRES-DR(200,150) and is much less expensive, because on the fine grid GMRES(100) has less orthogonalization.  However, Two-grid Deflated IDR is the best method.  It uses IDR(20)-Proj(150) on the fine grid.  It converges faster in terms of matrix-vector products and also is much less expensive per matrix-vector product.

\begin{figure}
\vspace{-2.65in}
\hspace{-.7in}
\includegraphics[scale=.75]{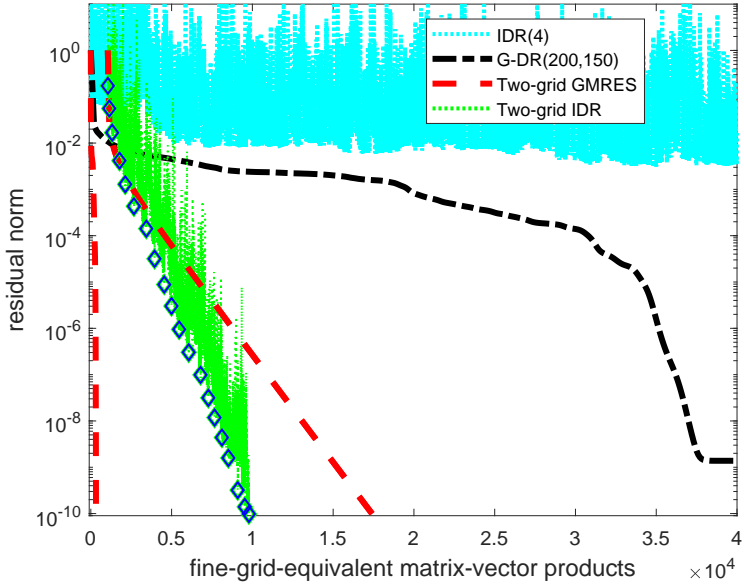}
\vspace{-2.8in}
\caption{ Matrix is from the simple Helmholtz equation with $\kappa = 100$.  Two-grid deflated GMRES uses GMRES-DR(200,150) on the coarse grid and GMRES(100)-Proj(150) on the fine grid.  Two-grid deflated IDR uses GMRES-DR(200,150) on the coarse grid and IDR(20)-Proj(150) on the fine grid; diamonds show the residual norm at the end of each of the 20 cycles.  Also compared are IDR and GMRES-DR(200,150).}
\end{figure}

We did try only deflating 100 eigenvalues as in Examples 1 and 2.  However, deflated IDR converges more than twice as slow.  For this difficult indefinite problem, many eigenvalues are needed for effective deflation.

Helmholtz problems are fairly complicated.  For example, if the wave number is increased and the fine grid uses the same discretization, then the coarse grid may need to be finer than in Example 4 in order to have good enough eigenvector approximations.  Our goal here is merely to show potential for the new approach; much more work is needed on Helmholtz problems for a thorough study.

{\it Example 5.}  For this example, the problem is the same as in the previous one, but we now use multigrid preconditioning.  Since multigrid does not converge for this matrix, the preconditioner has solution of linear equations from the positive shifted Laplacian with operator $- u_{xx} - u_{yy} + 100^2 u $~\cite{LaGi}.  With this positive shift, multigrid easily converges and thus can precondition the Helmholtz matrix (the negatively shifted Laplacian).  This preconditioning is used on both the fine grid and the coarse grid and makes the problem easier to solve; see Figure 5.2.  The problem is still indefinite and while BiCGStab still does not converge, now IDR does converge.  GMRES-DR(150,100) converges in less than half the number of matrix+preconditioner applications than IDR, but has more orthogonalization expense.
The two-grid deflated methods use GMRES-DR(150,100) on the coarse grid, stopping when 80 eigenpairs have converged to residual norm $10^{-8}$.  Then on the fine grid, deflated GMRES uses GMRES(100)-Proj(100) and deflated IDR uses IDR(5)-Proj(100) (the IDR converges almost 25\% slower with $ncyc = 20$ instead of 5).  Both deflated methods converge faster than the other methods, but the deflated IDR has much less orthogonalization.

\begin{figure}
\vspace{-2.65in}
\hspace{-.7in}
\includegraphics[scale=.75]{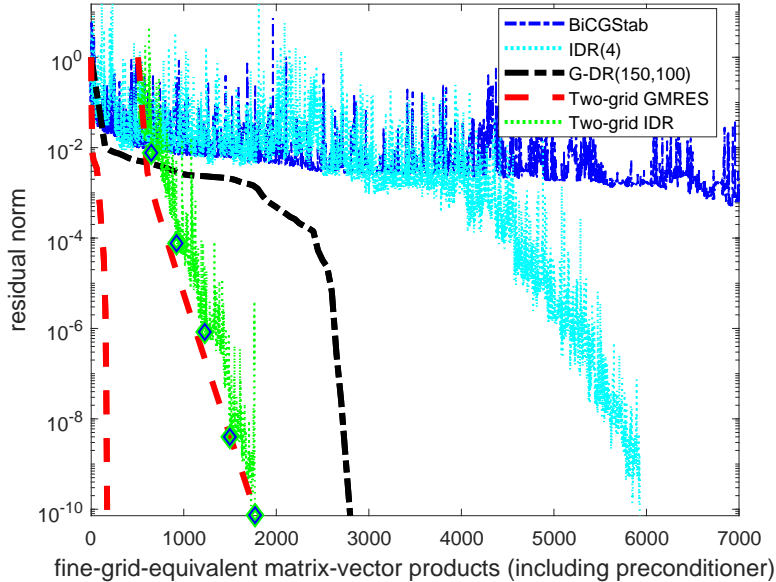}
\vspace{-2.8in}
\caption{Matrix is from the simple Helmholtz equation with $\kappa = 100$.  Now multigrid preconditioning is used.  Two-grid deflated GMRES uses GMRES-DR(150,100) on the coarse grid and GMRES(100)-Proj(100) on the fine grid.  Two-grid deflated IDR uses GMRES-DR(150,100) on the coarse grid and IDR(5)-Proj(100) on the fine grid; diamonds show the residual norm at the end of each of the 5 cycles.  These methods are compared to IDR and GMRES-DR(150,100). }
\end{figure}

In this experiment, the positively shifted Laplacian linear equations are solved accurately.  However, one can relax that and apply the multigrid only to partial convergence.  The methods still work, but the results vary some.

We also tried complex shifts for the multigrid preconditioner~\cite{ErOoVu06}.  There is faster convergence, but greater cost per iteration due to subspaces becoming complex. For a complex Helmholtz problem, this would not be a disadvantage and should be a subject for further study.

\subsection{Biharmonic Matrix}

We next consider matrices from discretizing a biharmonic differential equation.
Matrices from this differential equation quickly become very ill-conditioned as the discretization size gets small.
The biharmonic examples demonstrate first that residual norms can jump up during the projection over approximate eigenvectors.  Then the second example has faster convergence if approximate eigenvectors are improved on the fine grid.

{\it Example 6.}
The partial differential equation is $ - u_{xxxx} - u_{yyyy} + 40 u_{xxx}  = f$ on the unit square with zero boundary conditions.
The discretization has $h = 1/256$, giving a matrix of size $n = 65{,}025$.
The right-hand side is chosen random normal.  Due to the ill-conditioning, all Biharmonic tests have residual tolerance for the linear equations of only $10^{-8}$.  BiCGStab and IDR do not converge and GMRES-DR is slow and expensive.  We only give results for deflated BiCGStab.  The top half of Figure 5.3 shows results with three choices of $m$ and $k$ for GMRES-DR(m,k) on the coarse grid.  GMRES-DR(100,50) finds 40 eigenpairs to residual norms below $10^{-8}$, GMRES-DR(150,100) gets 80 to that tolerance and GMRES(200,150) stops when 120 have converged.  Deflated BiCGStab is used on the fine grid with $ncyc = 20$ cycles.  The convergence is plotted against cost for both coarse and fine grid phases.  Here each matrix-vector product is counted as 13 vector operations (the number of non-zeros in most rows).  The cost for the coarse grid phase is greater when more eigenvalues are computed, but then the convergence is faster on the fine grid.  Using 100 approximate eigenvectors is best overall for this example.  If there were multiple right-hand sides, then 150 would be clearly better for subsequent right-hand sides because of the faster fine grid convergence.

The lower half of Figure 5.3 has a portion of the fine grid restarted BiCGStab with $k=150$, plotted against matrix-vector products.  The residual norm jumps up by a significant amount during each projection over the approximate eigenvectors.  For instance, it increases from $2.5*10^{-6}$ to $1.2*10^{-5}$ in between cycles 14 and 15 (at matrix-vector product 3415).  As a result, the second of the two convergence criteria in part 1.\thinspace b) of the algorithm is activated.  Figure 5.4 shows eigencomponents of residual vectors for a smaller version of this problem.  The matrix is size $n = 961$ and the coarse grid matrix is size $49$.  Here 10 eigenvectors are computed accurately on the coarse grid and moved to the fine grid.  The top of the figure shows all 961 eigencomponents during part of a deflated BiCGStab run on the fine grid.  The (red) circles are from the residual vector after the projection at the start of the fifth cycle.  Then the (black) squares are after BiCGStab has been applied.  Finally, the (blue) dots are after the next projection at the start of the sixth cycle, and they mostly overlie the circles but are a little better on average.  Most of the components increase dramatically with the projection.  Then, fortunately they are reduced by the Krylov iteration.  The lower part of the figure shows that some of the components corresponding to the small eigenvalues are reduced by the projection.  This reduction is important because they are for the most part not reduced by the BiCGStab.  This essential reduction of the small eigencomponents by the projection makes up for the increase of the other components, because overall this deflated method does better than regular non-restarted BiCGStab.

\begin{figure}
\vspace{-2.65in}
\hspace{-.7in}
\includegraphics[scale=.75]{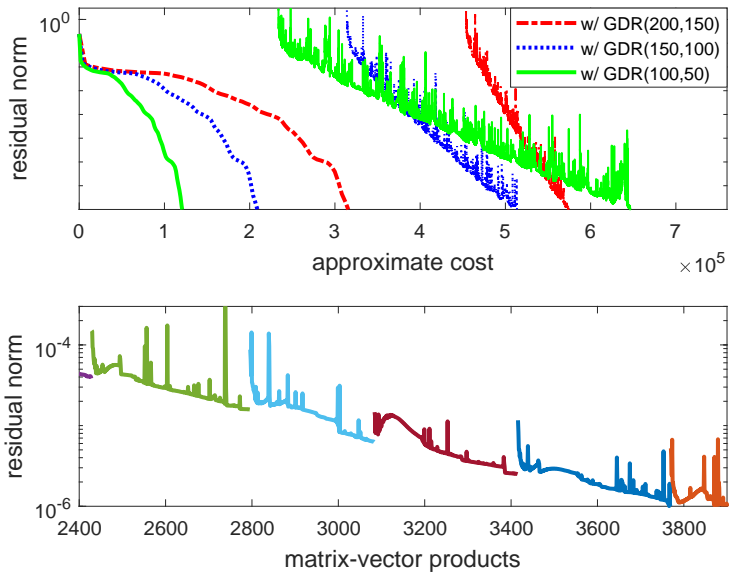}
\vspace{-2.8in}
\caption{ Matrix is from a biharmonic differential equation.  The top half has deflated BiCGStab, with convergence shown on both the coarse and fine grids.  The coarse grid has GMRES-DR(m,k) for three choices of m and k.  The fine grid has BiCGStab(50)-Proj(k) for the three values of k.  The bottom half has a view of part of the fine grid convergence with $k=150$, showing the jump in residual norm with each projection. }
\end{figure}

((think it is BiCGStab(20)-Proj(k) in Fig 5.3, not 50))

\begin{figure}
\vspace{-2.65in}
\hspace{-.7in}
\includegraphics[scale=.75]{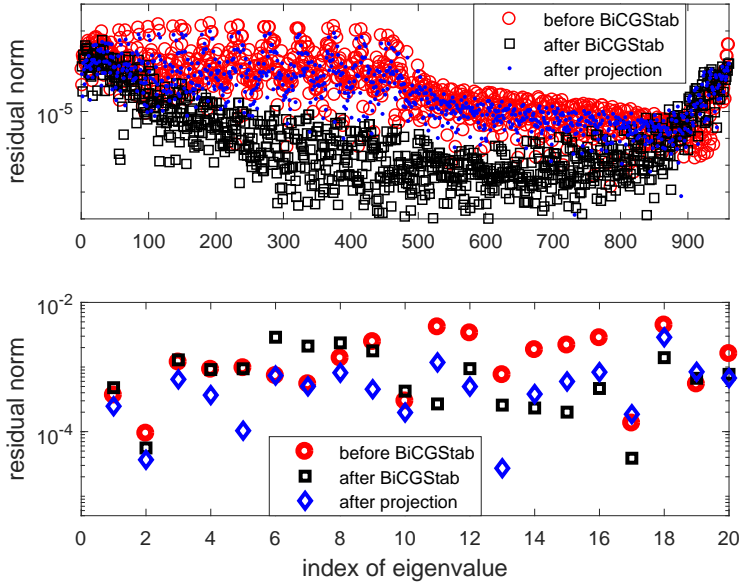}
\vspace{-2.8in}
\caption{ Small matrix from a biharmonic differential equation.  The top portion of the figure has eigencomponents of the residual vector before applying BiCGStab in the fifth cycle (red circles), then after the BiCGStab (black squares), and then after the projection over the approximate eigenvectors (blue dots).  The bottom portion has only the first 20 eigencomponents. }
\end{figure}

\subsection{Preconditioned Biharmonic}

We continue the biharmonic example, but now use an incomplete LU factorization.  This example also considers improving the eigenvectors on the fine grid (phase 2 of the Two-grid Deflated GMRES algorithm).

{\it Example 7.}
We use incomplete factorization preconditioning for the biharmonic differential equation from the previous example.  This is done with Matlab's ``ilu" command with no fill, after adding 0.5 to all diagonal elements of the matrix (this is needed to make the preconditioning effective).  We choose a smaller discretization size than the previous example, because the preconditioning allows us to solve a harder problem.  The matrix has $n = 511^2 = 261{,}121$ and the course grid is size $127^2 = 16,129$.  An ILU preconditioner is generated for both the fine and course grids.  As in the previous example, $rtol$ is set at $10^{-8}$.

Table 5.1 has results of a few tests.  The second row of the table (``cycles phase 1") has the number of GMRES-DR(150,100) cycles on the coarse grid.  The third row has the cycles of Arnoldi-E on the fine grid used to improve the approximate eigenvectors (phase 2.\thinspace of the Two-grid Deflated GMRES algorithm).  The fourth row has the number of cycles on the fine grid; this happens to be exactly 200 for GMRES-Proj and is specified to be 50 for restarted BiCGStab.  The fifth row of the table has the number of applications of matrix-vector product plus preconditioner on the fine grid.  Then the next row adds to this the coarse grid mvp's+preconditionings, scaled by 16 to give the fine-grid-equivalent total.  The last row has the approximate cost which counts 26 for each matrix-vector product plus preconditioner (13 for the mvp and 13 for preconditioning) and then adds the vector ops.  The coarse grid costs are included, but again scaled by 16.
The first result in the table (the second column) is for regular, non-restarted BiCGStab.  For convergence, this takes 34,907 applications of a matrix-vector product plus a preconditioning.  The next column has Two-grid Deflated GMRES, with GMRES-DR(150,100) on the coarse grid and GMRES(50)-Proj(100) on the fine grid.  This is much better than BiCGStab in terms of the numbers of applications of matrix-vector product plus preconditioner, using the fine-grid-equivalent if 10,203 of them.  However, the cost estimate in the last row of the table is only a little better (854 thousand compared to 1152 thousand) due to the orthogonalization in the 200 cycles of GMRES(50).  If GMRES(100)-Proj(100) is used instead, the number of cycles on the fine grid reduces to 73, and thus the fine-grid equivalent matrix-vector product plus preconditioners goes down to 7576.  However the cost goes up to 988 thousands.  The next column has Deflated BiCGStab.  The same GMRES-DR run is used on the coarse grid and then the fine grid has BiCGStab(50)-Proj(100), with 50 cycles of BiCGStab and projections over the 100 approximate eigenvectors in between.  The cost is significantly reduced compared to the other methods.  We note that with $ncyc = 20$, there are 12,168 mvp+prec's.  For this example, more restarts are needed, perhaps because the approximate eigenvectors are not as good and so need to be deflated more often.  The next to last column has Arnoldi-E(150,100) applied to the approximate eigenvectors for 10 cycles with the 10 smallest approximate eigenvectors used as starting vectors.  Then the last column has the Arnoldi-E applied as in~\cite{MGArn} until the first 80 approximate eigenvectors have residual norms below $10^{-3}$ (before this improvement, the 80 have residual norms from $3.3*10^{-4}$ to $2.4*10^{-2}$).  This improvement takes 41 cycles.  Both of these raise the cost but reduce the BiCGStab iterations.  For cases with multiple right-hand sides, this eigenvector improvement could be worthwhile.  Improving eigenvectors was generally not beneficial in the earlier examples.  Probably the best way to decide if eigenvectors should be improved is by experimenting.

\begin{table}

\caption{Results for biharmonic matrix preconditioned by ILU factorization.  The first method is regular non-restarted BiCGStab.  Second is Two-grid Deflated GMRES with GMRES-DR(150,100) on the coarse grid and GMRES(50)-Proj(100) on the fine grid.  The last three columns have the same coarse grid work and then restarted BiCGStab(50)-Proj(100) on the fine grid.  The first of these three tests has no improvement of the approximate eigenvectors on the fine grid, and the last two have 10 and 41 cycles of improvement respectively. }

\begin{center}
\begin{tabular}{|c|c|c|c|c|c|}  \hline\hline
fine grid method       & BiCGSt    & G-Proj    & Defl. Bi. & Defl. Bi. & Defl. Bi. \\ \hline \hline
cycles phase 1  & -         & 31        & 31        & 31        & 31        \\ \hline
cycles phase 2  & -         & 0         & 0         & 10        & 41        \\ \hline
cycles phase 3  & 1         & 200       & 50        & 50        & 50        \\ \hline
f. g. mvp + prec & 34,907   & 10,000    & 7877      & 6074      & 3069      \\ \hline
total f. g. equiv. mvp+ & 34,907 & 10,203    & 8080      & 6827      & 5372      \\ \hline
cost (in thousands)& 1152      & 854       & 322       & 801       & 2369      \\ \hline
 \hline

\end{tabular}
\end{center}
%\label{}

\end{table}

\subsection{QCD Problem}

An important computational area that especially needs further development of linear solvers is quantum chromodynamics (QCD).  Increasingly large and difficult systems of equations are developed in QCD and so the methods need to increase in effectiveness.  Multigrid methods have been developed for QCD~\cite{Babich:2010qb,Cohen:2012sh,Brannick:2014vda}.  They show potential, but are not as efficient as multigrid for standard PDE's and are much more complicated to implement.  Currently multigrid methods are not standard for QCD but are one possible tool.  We wish to add Two-grid Deflated BiCGStab to the possible methods.  For this we will use some of the framework that has been previously developed for QCD multigrid.  Here we do preliminary testing with a moderately sized problem in the simpler QCD situation of the two-dimensional Schwinger model~\cite{Schwinger:1962tp}.  It is beyond the scope of this paper to test other models or to give a full comparison with other QCD methods.

{\it Example 8.}  As mentioned above, the matrix is from the Schwinger model.  It is size $n=294{,}912$.  A coarse grid matrix is developed of size $n=9216$ using techniques described in~\cite{WHYTEWM}.  Along with this, a prolongation operator is formed for moving vectors from the coarse to the fine grid.  The original fine grid matrix has 9 non-zeros per row, while the coarse grid matrix has 80.  There are significant costs to forming the coarse grid matrix, including 1128 matrix-vector products with the fine grid matrix to develop very rough approximations to small eigenvectors that are needed in the process of forming it.  There are of course also costs for GMRES-DR on the coarse grid.  However, these costs may not be significant if many right-hand sides are solved as is commmon for QCD problems.  Therefore, here we only compare matrix-vector products during the solve phase on the fine grid (Phase 3).  The matrix is shifted by three different values, 0.061, 0.062 and 0.063, where critical mass is roughly at 0.062.  Since this is a relatively small and easy problem, the last shift is slightly past critical in order to simulate a more difficult problem.

The first row of results in Table 5.2 is for regular BiCGStab, and it does not converge for the most difficult matrix.  It is not unusual for BiCGStab to fail for difficult QCD problems.  Then three tests of Two-grid Deflated BiCGStab are given with increasing work on developing the approximate eigenvectors.  The second row of results has GMRES-DR(80,40) for 20 cycles on the coarse grid and no Phase 2 improvement.  This is better than regular BiCGStab for the first shift of 0.061, slower for the second shift and also does not converge for the third.  The third row of results has an added 15 cycles of improvement on the fine grid, targeting only the smallest five eigenvalues and eigenvectors (this is because in testing, it is more important to have accurate eigenvectors for the eigenvalues near zero).  The cost for this Phase 2 work is significant for one right-hand side, but is not as important if many are solved.  But also, we are now able to get convergence for the difficult shift.  The last row has the same method except there is more work on the coarse grid with 43 cycles (this is enough to have the five smallest eigenpairs converge to accuracy of $10^{-8}$).  These experiments show the new approach has potential for difficult problems.

\begin{table}

\caption{Solving a 2-D QCD matrix.  GMRES(80,40) on the coarse grid.  Fine grid matrix-vector products are given for solving one right-hand side, not including set-up costs for forming the coarse grid matrix.}

\begin{center}
\begin{tabular}{|c|c|c|c|c|}  \hline\hline
\backslashbox{Method}{Shift}    & 0.061     & 0.062     &  0.063  \\  \hline \hline
% Method                                                     \\ \hline
non-restarted BiCGStab                      & 4880  & 6187      &   -      \\ \hline
Defl. Bi., 20 cyc phase 1, 0 cyc phase 2    & 3589  & 12,284    &   -   \\ \hline
Defl. Bi., 20 cyc phase 1, 15 cyc phase 2   & 3247  & 4002      & 4519  \\ \hline
Defl. Bi., 43 cyc phase 1, 15 cyc phase 2   & 2725  & 3354      & 3055  \\ \hline
 \hline

\end{tabular}
\end{center}
%\label{}

\end{table}

\section{Conclusion}

We have proposed a two-grid method that finds approximate eigenvectors with the coarse grid and uses them to deflate eigenvalues for linear equations on the fine grid.  This includes deflation for BiCGStab and IDR using only approximate right eigenvectors and novel use of restarting for these normally non-restarted nonsymmetric Lanczos methods.  This two-grid deflation is a very efficient way to deflate eigenvalues, because the difficult work of finding approximate eigenvectors is done for an easier problem on the coarse grid.  This is particularly useful for multiple right-hand side problems, because the coarse grid work only needs to be done once and then can be applied for all right-hand sides.

The new two-grid deflated methods are not as efficient as multigrid when it works well, but the new methods can potentially be used in situations where multigrid is not effective.  Additionally, they can be combined with multigrid preconditioning.

Many facets of these two-grid deflated methods could use further investigation.  For instance, three dimensional problems may have greater potential because the coarse grid matrix can be relatively smaller compared to the fine grid matrix.  QCD problems could use much further investigation, including going to four-dimensional problems.  Other possible future work is to use more grid levels; see~\cite{MGArn} for a multiple grid method for computing eigenvalues.

\bibliography{MGLinEq2}

\end{document}